\documentclass{aptpub}

\authornames{Torrey Johnson and Edward C Waymire} 
\shorttitle{Tree Polymers in the Infinite Volume Limit at Critical
Strong Disorder} 

\usepackage{MnSymbol}
\usepackage{url}
\usepackage{graphicx}

\numberwithin{equation}{section}

\numberwithin{equation}{section}  

\begin{document}

\title{Tree polymers in the  infinite volume limit at \\ critical
strong disorder} 

\authorone[Oregon State University]{Torrey Johnson}
\addressone{Department of Mathematics, Oregon State Univeristy, Corvallis, Oregon 97331}
\emailone{johnsotor@science.oregonstate.edu}
\authortwo[Oregon State Univeristy]{Edward C Waymire}
\addresstwo{Department of Mathematics, Oregon State Univeristy, Corvallis, Oregon 97331}
\emailtwo{waymire@math.oregonstate.edu}

\begin{abstract}
The a.s. existence of a polymer probability in the infinite volume limit  is readily obtained under
general conditions of weak disorder from standard theory on multiplicative cascades or 
branching random walk.   However, speculations in the case of strong disorder have
been mixed.   In this note existence of an infinite volume probability is established at critical strong disorder for which one has convergence in probability.   Some
calculations in support of a specific formula for the a.s.  asymptotic variance of the polymer path 
under strong disorder are also provided.
\end{abstract}



\section{Introduction and Preliminaries}
Polymers are abstractions of chains of molecules
embedded in a solvent by non-self-intersecting
polygonal paths of points whose probabilities are
themselves random (reflecting impurities of the 
solvent).   In this connection,
tree polymers take advantage of a
particular way to determine path structure and
their probabilities as follows.

Three different references to paths occur in this formulation.
An \emph{$\infty$-tree path} is a sequence $s=(s_1,s_2,\ldots) \in \{-1,1\}^\mathbb{N}$ emanating from a root 
$0$.
A \emph{finite tree path} or \emph{vertex} $v$ is a finite sequence $v=s|n=(s_1,\ldots,s_n)$, read 
\lq\lq path $s$ restricted
to level $n$\rq\rq,
of length $|v|=n$.  The symbol $*$ denotes concatenation of finite tree paths; if $v=(v_1,
\ldots,v_n)$ and $t=(t_1,\ldots,t_m)$, then $v*t=(v_1,\ldots,v_n,t_1,\ldots,t_m)$.
Vertices belong to $T:=\bigcup_{n=0}^\infty \{-1,1\}^n$, and can be viewed as unique finite
paths to the root of the directed binary tree $T$
equipped  with the obvious graph structure.
We also write
\begin{equation*}
\partial T = \{-1,1\}^\mathbb{N}
\end{equation*}
for the boundary of $T$.  The third type of path, and the one of main interest to polymer
questions, is that of the
 \emph{polygonal tree path} defined by $n\rightarrow(s)_n:=\sum_{j=1}^ns_j$, $n\geq0$, with
$(s)_0:=0$, for a given $s \in \partial T$.

$\partial T$ is a compact, topological Abelian group for coordinate-wise multiplication
and the product topology.  The \emph{uniform distribution} on $\infty$-tree paths is the
Haar measure on $(\partial T, \mathcal{B})$, i.e.
\begin{equation*}
\lambda(ds)=\left( \frac{1}{2}\delta_+(ds) + \frac{1}{2}\delta_-(ds) \right)^\mathbb{N}.
\end{equation*}
 
Let $\{X_v:v\in T\}$ be an i.i.d. family of positive random variables on $(\Omega,\mathcal{F},
P)$ with $\mathbb{E}X < \infty$; we denote a generic random variable with the common
distribution of $X_v$ by $X$.  Without loss of generality we may assume that $\mathbb{E}X=1$.
Define a sequence of \emph{random probability measures} $\text{prob}_n(ds)$ on $(\partial T, 
\mathcal{B})$ by the prescription that 
\begin{equation*}
\text{prob}_n(ds) << \lambda(ds)
\end{equation*}
with
\begin{equation*}
\frac{d\text{prob}_n}{d\lambda}(s)=Z_n^{-1}\prod_{j=1}^n X_{s|j}
\end{equation*}
where
\begin{equation*}
Z_n=\int_{\partial T}\prod_{j=1}^n X_{s|j}\lambda(ds)=\sum_{|s|=n}\prod_{j=1}^n X_{s|j}2^{-n}.
\end{equation*}

Observing that $\{Z_n:n=1,2\ldots\}$ is a positive martingale, it follows that
\begin{equation*}
Z_\infty := \lim_{n\rightarrow\infty}Z_n
\end{equation*}
exists a.s. in $(\Omega,\mathcal{F},P)$.  According to a classic theorem of Kahane and
Peyri\`ere (1976) in the context of multiplicative cascades, and Biggins (1976) in the 
context of branching random walks, one has the following dichotomy:
\begin{align*}
P(Z_\infty>0) = 1 \quad &\Longleftrightarrow \quad \mathbb{E}X \ln X < \ln 2   \\
P(Z_\infty=0) = 1 \quad &\Longleftrightarrow \quad \mathbb{E}X \ln X \geq \ln 2.  
\end{align*}
The a.s. occurance of the event $[Z_\infty > 0]$ is refered to as \emph{weak disorder} and
that of $[Z_\infty=0]$ as \emph{strong disorder}; see Bolthausen (1989).  In particular, 
the critical case $\mathbb{E}X\ln X =\ln 2$ is strong disorder.  In the case of tree polymers
one may view the notions
of weak/strong in terms of a disorder parameter
defined by
$\mathbb{E}X\ln X$ and relative to the branching 
rate, $\ln 2$.

In this short communication we provide some new insights
into a few delicate problems for the case of strong
disorder. 

\section{Tree Polymers under Weak Disorder}
To set the stage for contrast, we record a rather robust
consequence of weak disorder.
\begin{theorem}
Under weak disorder, there is a random probability measure $\text{\emph{prob}}_\infty(ds)$ on
$(\partial T ,\mathcal{B})$ such that a.s.
\begin{equation*}
\text{\emph{prob}}_n(ds) \Rightarrow \text{\emph{prob}}_\infty(ds)
\end{equation*}
where $\Rightarrow$ denotes weak convergence.
\end{theorem}

\begin{proof}
Define $\lambda_n(ds)=Z_n\text{prob}_n(ds)$, $n=1,2,\ldots$.  By Kahane's $T$-martingale theory, 
e.g.,  Kahane (1989), $\lambda_n(ds)$ converges vaguely to a non-zero 
random measure $\lambda_\infty(ds)$ on $(\partial T, \mathcal{B})$ with probability one.  By definition of weak disorder $Z_n \rightarrow Z_\infty > 0$ a.s., thus we obtain
\begin{equation*}
\text{prob}_n(ds)=Z_n^{-1}\lambda(ds) \Rightarrow Z_\infty^{-1}\lambda_\infty(ds) \quad \text{a.s.}
\end{equation*}
\end{proof}

Notice that in the case of no disorder, i.e. $X=1$ a.s., one has
\begin{equation*}
\text{prob}_n(ds)=\lambda(ds) \quad \forall n=1,2,\ldots.
\end{equation*}
Moreover, under $\lambda(ds)$, the polygonal paths are simply symmetric simple random walk paths,
where the probability theory is quite will-known and complete.  For example, the central limit 
theorem takes the form
\begin{equation*}
\lim_{n\rightarrow\infty} \lambda\left(\left\{s\in\partial T: \frac{(s)_n}{\sqrt{n}}\leq x
\right\}\right) = \frac{1}{\sqrt{2\pi}}\int_{-\infty}^x e^{-\xi^2/2}d\xi.
\end{equation*}
For probability laws involving convergence in distribution, one may ask if the CLT continues
to hold a.s. with $\lambda(ds)$ replaced by $\text{prob}_n(ds)$.  This form of universality
 was answered in the 
affirmative by Waymire and Williams (2010) for weak disorder under the additional assumption
that $\mathbb{E}X^{1+\delta}<\infty$ for some $\delta>0$.  
Problems involving limit laws  
such as a.s. strong laws, a.s. laws of the iterated logarithm, etc, however,
require an infinite volume probability $\text{prob}_\infty(ds)$ for their formulation.  While the 
preceding proposition answers this in the case of weak disorder,  the problem is open for
strong disorder.  Moreover, it has been speculated by Yuval Peres (private communication) that
$\text{prob}_n(ds)$ will a.s. have infinitely many weak limit points under strong disorder.  However, in the case of critical strong disorder we show that a natural infinite volume polymer exists and is related to the finite volume
polymers through limits in probability.

\section{Tree Polymers at Critical Strong Disorder}

In this section we show the  existence under critical strong disorder, i.e., assuming $\mathbb{E}X\ln X = \ln 2$,
 of an infinite volume polymer probability
 $\text{prob}_\infty(ds)$ that may be viewed as the weak
 limit in probability of the sequence $\text{prob}_n(ds), n\ge 1,$ in the sense that its characteristic function
 is the limit in probability of the corresponding 
 sequence of characteristic functions of 
 $\text{prob}_n(ds), n\ge 1$. 
 
For $v\in T$, $v=(v_1,\ldots,v_m)$, say, let
\begin{equation*}
\Delta_m(v)=\{s\in\partial T:s_i=v_i, \text{ } i=1,\ldots,m\}, \qquad |v|=m.
\end{equation*}
Since $T$ is countable there are countably many such finite-dimensional rectangles in $\partial T$.

For $m>n$, note that
\begin{align*}
\text{prob}_n(\Delta_m(v)) &= \int_{\Delta_m(v)}\frac{d\text{prob}_n}{d\lambda}(s)\lambda(ds) \\
                           &= \int_{\Delta_m(v)}Z_n^{-1}\prod_{j=1}^nX_{s|j}\lambda(ds) \nonumber\\
                           &= Z_n^{-1}\int_{\Delta_m(v)}\prod_{j=1}^nX_{v|j}\lambda(ds) \nonumber\\
                           &= Z_n^{-1}\prod_{j=1}^n X_{v|j}\cdot 2^{-m}.
\nonumber\end{align*}

For example,
\begin{align*}
\text{prob}_1(\Delta_m(v)) &= Z_1^{-1}X_{v|1}2^{-m}, \qquad Z_1=\frac{X_+ + X_-}{2} \\
													 &= \frac{X_{v|1}2^{-(m-1)}}{X_+ + X_-}\\
													 &= \left\{
													    \begin{array}{rl}
													    \frac{X_+2^{-(m-1)}}{X_+ + X_-}, & \quad v|1=+1 \\
													    \frac{X_-2^{-(m-1)}}{X_+ + X_-}, & \quad v|1=-1.
													    \end{array}\right.
\end{align*}
$\sum_{|v|=m}\text{prob}_1(\Delta_m(v))=1$ since there are $2^m$ such $v$'s, half of which
have $v_1=+1$ and the other half have $v_1=-1$.

For $m\leq n, |v|=m$, we have
\begin{align*}
\text{prob}_n(\Delta_m(v)) 
&= Z_n^{-1}\int_{\Delta_m(v)}\prod_{j=1}^n X_{s|j} \lambda(ds) \\
&= Z_n^{-1}\prod_{j=1}^m X_{v|j}\sum_{|t|=n-m}\prod_{j=1}^{n-m} X_{(v*t)|j}2^{-n}\nonumber\\
&= Z_n^{-1}\left(\prod_{j=1}^m X_{v|j} 2^{-m}\right)Z_{n-m}(v),
\nonumber\end{align*}
where
\begin{equation*}
Z_0(v)=1, \quad Z_{n-m}(v)=\sum_{|t|=n-m}\prod_{j=1}^{n-m}X_{(v*t)|j}2^{-(n-m)}.
\end{equation*}
In particular, $Z_n=Z_n(0)$, where $0\in T$ is the root.

Note that
\begin{align*}
Z_n &= \sum_{|u|=m}\sum_{|t|=n-m}\prod_{j=1}^{m}X_{u|j}2^{-m}\prod_{j=1}^{n-m}X_{(u*t)|j}2^{-(n-m)} \\
	  &= \sum_{|u|=m}Z_{n-m}(u)\prod_{j=1}^m X_{u|j}2^{-m}.
\nonumber\end{align*}
Thus, letting $a_k = 1/\sqrt{k}, k\ge 1$,
\begin{align*}
\text{prob}_n(\Delta_m(v))
&= \frac{D_{n-m}(v)\prod_{j=1}^m X_{v|j}2^{-m}\frac{Z_{n-m}(v)}{a_{n-m}D_{n-m}(v)}}
{\sum_{|u|=m}D_{n-m}(u)\left(\prod_{j=1}^m X_{v|j}2^{-m}\right)\frac{Z_{n-m}(u)}{a_{n-m}D_{n-m}(u)}}\\
&\longrightarrow \frac{D_\infty(v)\prod_{j=1}^m X_{v|j}2^{-m}}
{\sum_{|u|=m}D_\infty(u)\left(\prod_{j=1}^m X_{v|j}2^{-m}\right)}
\nonumber\end{align*}
where (i) the convergence to $D_\infty(v)$ is the
almost sure 
limit of the {\it derivative martingale} obtained by
Biggins and Kyprianou (2004),
and (ii) $\lim_{n\rightarrow\infty}\frac{Z_{n-m}(v)}
{a_{n-m}D_{n-m}(v)} = c > 0$ is the limit in probability
at critical strong disorder recently obtained by Aid\'{e}kon
and Shi (2011).  The  constant 
$c = ({2\over\pi\sigma^2})^{1/2}$, for $\sigma^2 = 
\mathbb{E}\{X(\ln(X))^2\}  - (\mathbb{E}\{X\ln(X)\})^2 > 0$,
does not
depend on $v\in T$.   Aid\'ekon
and Shi (2011) also point out that the almost
sure positivity of  $D_\infty(v)$ follows from 
Biggins and Kyprianou (2004) and Aid\'ekon (2011) 
The sequence $a_k = k^{-{1\over 2}}, k\ge 1,$  is
referred to as the Seneta-Heyde scaling.  

\begin{remark} For each $v\in T$, there is a set $N(v)$ of probability zero such that
\begin{equation*}
D_\infty(v,\omega) = \lim_{n\rightarrow\infty}D_n(v,\omega), \quad \omega\in\Omega
\backslash N(v).
\end{equation*}
Since $T$ is countable, the set $N=\bigcup_{v\in T}N(v)$ is still a $P$-null subset of $\Omega$.  The almost sure
convergence of the derivative martingales is essential
to the construction of $\text{\emph{prob}}_\infty$ given
in the lemma below.
\end{remark}

We now define
\begin{equation*}
\text{prob}_\infty(\Delta_m(v),\omega)
=\frac{D_\infty(v,\omega)\prod_{j=1}^m X_{v|j}(\omega)2^{-m}}
{\sum_{|u|=m}D_\infty(u,\omega)\left(\prod_{j=1}^m X_{u|j}(\omega)2^{-m}\right)}
\end{equation*}
for $\omega\in\Omega\backslash N$.

\begin{lemma}
$\text{\emph{prob}}_\infty(\Delta_m(v),\omega)$ extends to a unique probability on
$(\partial T,\mathcal{B})$ for each $\omega\in\Omega\backslash N$.
\end{lemma}

\begin{proof}
We use Caratheodory extension, taking careful advantage of the fact that the sets $\Delta(v)$,
$v\in T$, are both open and closed subsets of the compact set $\partial T$.  For $\omega\in
\Omega\backslash N$, $\text{prob}_\infty(\cdot,\omega)$ extends to the algebra generated by
$\{\Delta(v):v\in T\}$ by addition.  Since $\partial T$ is compact and the rectangles are both
open and closed, countable additivity on this algebra must hold as a consequence of finite
additivity; i.e. if $\bigcup_{i=1}^\infty \Delta(v_i)$ is contained in the algebra
generated by $\{\Delta(v):v\in T\}$, then $\bigcup_{i=1}^\infty \Delta(v_i)$ is closed,
hence compact, and its own open cover, i.e. $\bigcup_{i=1}^\infty \Delta(v_i)
=\bigcup_{i=1}^l \Delta(v_{i_l})$ for some finite subsequence $\{i_j\}_{j=1}^l$ of
$\{1,2,\ldots\}$.
\end{proof}

\begin{theorem}
At critical strong disorder, for each finite
set $F\subseteq\mathbb{N}$
\begin{equation*}
\widehat{\text{\emph{prob}}_n(F)} \quad \Rightarrow \quad \widehat{\text{\emph{prob}}_\infty(F)} \qquad \text{in probability},
\end{equation*}
where $\widehat{\text{\emph{prob}}}_n, n\ge 1,
\widehat{\text{\emph{prob}}}_\infty$ denote their respective Fourier
transforms as 
probabilities on the compact abelian multiplicative
group $\partial T$ for the product topology.
\end{theorem}

\begin{proof}
The continuous characters of the group $\partial T$ are given by
\begin{equation*}
\chi_F(t) = \prod_{j\in F}t_j \quad \text{for finite sets } F\subseteq\mathbb{N}.
\end{equation*}
In particular there are only countably many characters of $\partial T$.  From standard
Fourier analysis it follows that we need only show that
\begin{equation*}
\lim_{n\rightarrow\infty}\mathbb{E}_{\text{prob}_n}\chi_F = \mathbb{E}_{\text{prob}_\infty}\chi_F
\quad \text{in probability}
\end{equation*}
for each finite set $F\subseteq\mathbb{N}$.  Let $m=\text{max}\{k:k\in F\}$.  Then for $n>m$,
\begin{align*}
\mathbb{E}_{\text{prob}_n}\chi_F &= \int_{\partial T=\bigcupdot_{|v|=m}\Delta_m(v)}\chi_F(s)\frac{d\text{prob}_n}{d\lambda}(s)\lambda(ds)\\
&= \sum_{|v|=m}\left(\prod_{j\in F}v_j\right) Z_n^{-1}(0) \prod_{j=1}^m X_{v|j}2^{-m}\sum_{|t|=n-m}
\prod_{j=1}^{n-m} X_{(v*t)|j}2^{-(n-m)}\\
&= \sum_{|v|=m}\left(\prod_{j\in F}v_j\right)\prod_{j=1}^m X_{v|j}2^{-m}
\frac{Z_{n-m}(v)}{Z_n(0)}\\
&= \sum_{|v|=m}\left(\prod_{j\in F}v_j\right)\prod_{j=1}^m X_{v|j}2^{-m}D_{n-m}(v)
\frac{\frac{Z_{n-m}(v)}{a_{n-m}D_{n-m}(v)}}{\sum_{|u|=m}\prod_{j=1}^m X_{u|j}2^{-m}D_{n-m}(u)
\frac{Z_{n-m}(u)}{a_{n-m}D_{n-m}(u)}}\\
&\longrightarrow \mathbb{E}_{\text{prob}_\infty}\chi_F,
\end{align*}
where the convergence is almost sure for terms of 
the form $D_{n-m}$ 
and in probability for those of the form
$Z_{n-m}/(a_{n-m}D_{n-m})$
as $n\to\infty$.
\end{proof}

\section{Diffusivity Problems at Strong Disorder}
 With regard to the aforementioned a.s. limits in distribution of polygonal tree paths, 
 Waymire and Williams (2010) also obtained a.s. 
limits of the form 
\begin{equation*}
\lim_{n\to\infty}{\ln E_{prob_n}e^{r(S)_n}\over n} = F(r)
\end{equation*}
under both weak and strong disorder.   Let us refer to these as almost sure 
{\it Laplace rates} in reference to the Laplace principle of large deviation theory.

In the case of weak disorder the universal
limit is $F(r) = \ln\cosh(r)$,  in a neighborhood of the origin,
 otherwise independent of the distribution
of $X$.  In addition to being independent of the distribution of $X$ within the range of
weak disorder, this universality of Laplace rates 
is manifested in the coincidence with the same limit
obtained for $X \equiv 1$, i.e., for simple symmetric random walk.  

 For an illustrative case of strong disorder, consider $X = e^{\beta Z - {\beta^2\over 2}}$, where
$Z$ is standard normal and $\beta \ge \beta_c = \sqrt{2\ln2}$. Then
from Waymire and Williams (2010), it follows that  a.s. in a neighborhood
of the origin that
\begin{equation*}
F(r) = r\tanh(rh(r)) + \beta^2h(r) - \beta\beta_c,
\end{equation*}
where $h(r)$ is the uniquely determined solution to
\begin{equation*}
\beta^2 h^2(r) + 2rh(r)\tanh(rh(r))
-2\ln\cosh(rh(r)) = \beta_c^2;
\end{equation*}
also see Waymire and Williams (Sec 6, Cor 2, 2010) for the general formulae in the case of 
strong disorder.    In particular, the universality of the Laplace rates breaks down, even
at critical strong disorder. 
A graph of $F(r)$ computed from MATLAB is 
indicated in Figure 1 for the strong disorder case of
$\beta = 2\beta_c$.  

\begin{figure}[ht]
\begin{center}
\includegraphics[scale=0.60]{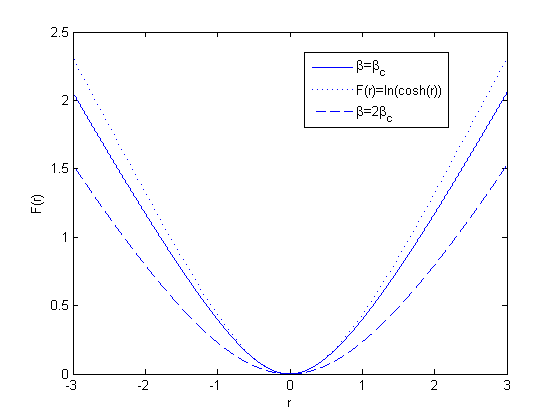}
\caption{Graph of the function $F$ for various $\beta$.}\label{Figure 1: }
\end{center}
\end{figure}

 Using the  equations defining $F(r)$ 
one may easily verify that $F(0) = 0, F^\prime(0) = 0$
and $F^{\prime\prime}(0) = { 2\beta\beta_c-\beta_c^2\over\beta^2}$.  
While these specific calculations follow directly from the general results of
 Waymire and Williams (2010),
from here one is naturally lead to speculate \footnote{To
 avoid potential confusion, let us mention that other forms of polymer scalings appear in the recent probability literature under which the polymer is referred to as \lq\lq superdiffusive\rq\rq even in the context of weak disorder; e.g.,  in reference to  wandering exponents in Bezerra, Tindel, Viens (2008).}  
that the asymptotic variance under strong disorder is obtained 
under diffusive scaling by $\sqrt{n}$ precisely as
\begin{equation*}
\sigma^2(\beta) = {2\beta\beta_c-\beta_c^2\over\beta^2},\quad \beta\ge\beta_c.
\end{equation*}
In particular this formula  continuously extends the weak disorder variance
$\sigma^2(\beta) \equiv 1, \beta < \beta_c,$ across $\beta = \beta_c$.  In any case,
this quantity is a basic parameter of the rigorously proven limit $F(r)$.

\section{Acknowledgment}
The authors thank the referee for spotting a serious error in the
original draft and providing the reference to
Aid\'ekon and Shi (2011) used in this paper.   The first author was partially supported by
an NSF-IGERT-0333257 graduate training grant in ecosystems informatics at Oregon State University, and the second author was partially supported by
a grant DMS-1031251 from the National Science Foundation.

%
%
%
%

\end{document}